\documentclass{article}

\usepackage{rotating}

\usepackage{amsfonts,amssymb,amsmath}

\def\BAL#1\EAL{\begin{align}#1\end{align}}

\newtheorem{theo}{Theorem}
\newcommand{\ncd}{\newcommand}
\ncd{\Eq}[1]{Eq.~(\ref{#1})}

\DeclareMathOperator{\lcm}{lcm}
\DeclareMathOperator{\glcd}{glcd}

\DeclareMathOperator{\Nr}{N}
\def\R{\mathbb R}
\def\N{\mathbb N}

\def\Z{\mathbb Z}
\def\I{\mathbb I}
\def\J{\mathbb J}
\def\H{\mathbb H}
\ncd{\gvct}[1]{{\boldsymbol{#1}}}
\newcommand{\nn}{\nonumber}

\begin{document}

\markboth{Michael Baake and Peter Zeiner}{Coincidences in 4 dimensions}

\title{Coincidences in 4 dimensions}

\author{Michael Baake and Peter Zeiner\thanks{$^\ast$Corresponding
author. Email: pzeiner@math.uni-bielefeld.de}$^\ast$\\\vspace{6pt} 
Faculty of Mathematics, Bielefeld University, Box 100131,
	33501 Bielefeld, GERMANY
}

\maketitle

\begin{abstract}
The coincidence site lattices (CSLs) of prominent 4-dimensional
lattices are considered.
CSLs in 3 dimensions have been used for decades to describe grain
boundaries in crystals. Quasicrystals suggest to also look at CSLs in
dimensions $d>3$. Here, we discuss the CSLs of the 
root lattice $A_4$ and the hypercubic
lattices, which are of particular interest both from the mathematical
and the crystallographic viewpoint. Quaternion algebras are used to derive
their coincidence rotations and the CSLs. We make use of the fact that the CSLs
can be linked to certain ideals and compute their indices, their
multiplicities and encapsulate all this in generating functions in terms of
Dirichlet series. In addition, we sketch how these results can be generalised
for $4$--dimensional $\Z$--modules by discussing the icosian ring.
\end{abstract}

\section{Introduction}

Let us begin with a brief recapitulation of concepts and methods.

\subsection{Coincidence site lattices and modules}
Let
$\varGamma\subseteq\R^d$ be a $d$-dimensional lattice and
$R\in \mathrm{O}(d)$ an isometry. Then, $R$ is called a
\emph{coincidence isometry} of $\varGamma$
if $\varGamma(R):=\varGamma\cap R\varGamma$ is a
lattice of finite index in $\varGamma$, and
$\varGamma(R)$ is called an (ordinary or simple)
\emph{coincidence site lattice (\mbox{CSL})} (for an introduction,
see~\cite{baa97}).
The group of all coincidence isometries of $\varGamma$ is called
$OC(\varGamma)$, whereas the subgroup of all orientation preserving
isometries is called the group of coincidence rotations and referred to
as $SOC(\varGamma)$.
The \emph{coincidence index} $\varSigma(R)$ is defined as the index of
$\varGamma(R)$ in $\varGamma$, and
$\sigma=\{\varSigma(R)\mid R\in OC(\varGamma)\}$ denotes the
(ordinary or simple) coincidence spectrum. These definitions can be
generalised to include the possibility of multiple CSLs,
see~\cite{baagri05,pzmcsl2,pzmcsl3}.

All these concepts also extend to $\Z$--modules~\cite{baa97,plea96}. 
By a $\Z$--module of dimension $d$ and rank $k$,
we mean the $\Z$--span of $k$ rationally independent vectors in $\R^d$.

\subsection{Quaternions and rotations}

Rotations in $\R^4$ can be parameterised by a pair of 
quaternions~\cite{koecheng,val,vign}.
Recall
that the quaternion algebra $\H(\R)$ is the vector space $\R^4$ equipped with
a special (associative but non--commutative) product.
In standard notation, where $1,i,j,k$
form an orthonormal basis of $\R^4$, the product satisfies the defining
relations\footnote{Note that these relations are very similar to the relations
satisfied by the Pauli matrices. In fact, dividing the Pauli matrices
by the imaginary unit gives a representation of $i,j,k$ by
$2\times2$--matrices.}
$i^2=j^2=k^2=-1$ and $ijk=-kji=-1$.
This definition is extended to all quaternions by linearity. All non--zero
quaternions $q=(a,b,c,d)=a+bi+cj+dk$ have an inverse
$q^{-1}=\frac{1}{|q|^2}\bar{q}$, where
$\bar{q}=(a,-b,-c,-d)$ is the conjugate of $q$ and
$|q|^2=a^2+b^2+c^2+d^2$ is its (reduced) norm. With this notation at hand,
we can write any rotation in $\R^4$ as
\BAL
Rx=R(q,p)x=\frac{1}{|qp|}qx\bar{p}.
\EAL
Reflections and, more generally, all orientation reversing isometries
may be parameterised by a pair of quaternions as well, but a
conjugation is required in addition.
For our purposes, it is sufficient to focus the
discussion to (proper) rotations. This is no restriction,
since all the lattices and modules we
discuss have a reflection as a symmetry operation and thus every
orientation reversing coincidence isometry can be written as the product of
this particular reflection
and a coincidence rotation with the same index. In other words, every
orientation reversing coincidence isometry is symmetry related to
a coincidence rotation with the same index (see below).
Hence all CSLs are already
obtained by (proper) coincidence rotations and the spectrum does not change
either. There are exactly as many coincidence rotations as orientation
reversing isometries, and the total number of coincidence isometries
is just twice the number of coincidence rotations.

\subsection{Symmetry related rotations and CSLs}

Two coincidence rotations $R$ and $R'$ are called
\emph{symmetry related}
if there is an element $Q$ in the point group $P(\varGamma)$ such that
$R'=RQ$. Note that $R\varGamma=R'\varGamma$ holds if and only if
$R$ and $R'$ are symmetry related. In particular, two symmetry related
rotations $R$ and $R'$ generate the same CSL, so that
$\varGamma(R)=\varGamma(R')$. However, the converse statement is not
true in general, i.e., two rotations generating the same CSL need \emph{not} be
symmetry related. In particular, in all the examples presented below, there are
additional rotations that generate the same CSL.

\section{Hypercubic lattices}

The hypercubic lattices have been discussed in detail in~\cite{baa97,pzcsl2}.
Here, we just review the most important facts, reformulate some
of them and add some new details.
In particular, we present an explicit expression for the number of CSLs.

\subsection{Centred hypercubic lattice}

The centred hypercubic lattice $D_4^*$, which is the dual lattice of the 
root lattice $D_4$, which in turn is a similar sublattice of $D_4^*$, can be
identified with the Hurwitz ring $\J$ of integer quaternions. This
is the $\Z$--span (meaning the set of all integral linear combinations)
of the quaternions
$(1,0,0,0), (0,1,0,0), (0,0,1,0),  \frac{1}{2}(1,1,1,1)$.
Any coincidence rotation can be parameterised by a pair of integer quaternions.
More precisely, a rotation is a coincidence rotation of $D_4^*$ if and only if
it can be parameterised by an admissible pair of
quaternions~\cite{baa97,pzcsl2}. 
Here,
we call a quaternion $q\in\J$ \emph{primitive} if $\frac{1}{n}q\in\J$
with $n\in\N$ implies $n=1$. A primitive pair $(q,p)$ is called
\emph{admissible} if $|qp|\in\Z$. In addition, we call a primitive quaternion
$q\in\J$ \emph{reduced} if $|q|^2$ is odd (for a more general discussion of
reduced quaternions, see~\cite{baaplea06}).
Due to the unique (left or right) prime decomposition in $\J$, we can decompose
any primitive quaternion $q\in\J$ as $q_r s$, where $q_r$ is
reduced and $|s|^2$ is a power of $2$. Note that $q_r$ is unique
up to right multiplication by a unit\footnote{Recall that the units are those
quaternions $u\in\J$ for which $|u|^2=1$.} of $\J$.

For any admissible
pair $(q,p)$, we can define the factors
\BAL 
\alpha_q:=\sqrt{\frac{|p_r|^2}{\gcd(|q_r|^2, |p_r|^2)}}&&\mbox{and}&&
\alpha_{p}:=\sqrt{\frac{|q_r|^2}{\gcd(|q_r|^2, |p_r|^2)}}.\nn
\EAL
The pair $(q_\alpha,p_\alpha):=
(\alpha_q q_r,\alpha_p p_r)$ is called the \emph{reduced extension} of the pair
$(q,p)$. Note that two admissible pairs are symmetry related
if and only if their reduced extension pairs are equal (up to units).
The coincidence index of a rotation $R(q,p)$ can now be
written~\cite{baa97,pzcsl2} as
\BAL
\varSigma_{D_4}(R(q,p))
=\lcm(|q_r|^2,|p_r|^2)
=|q_\alpha|^2
=|p_\alpha|^2.
\EAL
Since any (positive) integer can be written as a sum of four squares, this
implies that the spectrum $\sigma$ is the set of all odd natural numbers,
$\sigma=\{1,3,5,7,\ldots\}$

The CSLs can be calculated explicitly~\cite{pzcsl5}. They read
\BAL
\J\cap \frac{q\J\bar{p}}{|qp|}
= q_\alpha\J +\J\bar{p}_\alpha.
\EAL

Since the point group of $D_4$ contains $576$ rotations, the number of
coincidence rotations of a given index $n$ can be written as 
$576 f_{D_4}^{rot}(n)$.
This gives a multiplicative arithmetic function, i.e.
\BAL
f_{D_4}^{rot}(mn)&=f_{D_4}^{rot}(m)f_{D_4}^{rot}(n) && 
\mbox{if $m,n$ are coprime,}
\EAL
which is completely determined by $f_{D_4}^{rot}(1)=1$, $f_{D_4}^{rot}(2^r)=0$
for $r\geq 1$ and
\BAL
f_{D_4}^{rot}(p^r)&=\frac{p+1}{p-1}p^{r-1}(p^{r+1}+p^{r-1}-2)
\quad \mbox{if $p$ is an odd prime, $r\geq 1$.}
\EAL
The multiplicativity of $f_{D_4}^{rot}(n)$ suggests the use of
a Dirichlet series
as a generating function for it:
\BAL
\Phi_{D_4}^{rot}(s)&=\sum_{n=1}^\infty \frac{f_{D_4}^{rot}(n)}{n^s}
=\prod_{p\ne 2}\frac{(1+p^{-s})(1+p^{1-s})}{(1-p^{1-s})(1-p^{2-s})}\nn\\
&=1+\frac{16}{3^s}+\frac{36}{5^s}+\frac{64}{7^s}+\frac{168}{9^s}
+\frac{144}{11^s}+\frac{196}{13^s}+\frac{576}{15^s}+\frac{324}{17^s}+\cdots.
\label{frotpcub}
\EAL
Similarly, one can calculate the number $f_{D_4}^{}(n)$ of different CSLs of a
given index $n$. Clearly, one has $f_{D_4}^{}(n)\leq f_{D_4}^{rot}(n)$, and
we do not have equality in general,
since two CLSs that are not symmetry related
may generate the same CSL. In fact, we have the following theorem which tells
us when two CSLs are equal~\cite{pzcsl5}.
\begin{theo}
\label{equalJ}
Let $(q_1,p_1)$ and $(q_2,p_2)$ be two primitive reduced admissible pairs.
Then,
\BAL
\J\cap \frac{q_1 \J \bar{p}_1}{|q_1p_1|}=\J\cap \frac{q_2 \J \bar{p}_2}{|q_2p_2|}
\EAL
holds if and only if
$|q_1p_1|=|q_2p_2|$, $\lcm(|q_1|^2,|p_1|^2)=\lcm(|q_2|^2,|p_2|^2)$,
$\glcd(q_1,|p_1q_1|)=\glcd(q_2,|p_2q_2|)$ and 
$\glcd(p_1,|p_1q_1|)=\glcd(p_2,|p_2q_2|)$ hold.
\end{theo}
Here, $\glcd$ denotes the greatest left common divisor
in $\J$, which is defined up to units of $\J$. We find that $f_{D_4}^{}(n)$ is
multiplicative, too, and is determined by $f_{D_4}^{}(1)=1$,
$f_{D_4}^{}(2^r)=0$ for $r\geq 1$ and
\BAL
f_{D_4}^{}(p^r)&=
\begin{cases}
\frac{(p+1)^2}{p^3-1}\left(p^{2r+1}+p^{2r-2}-2p^{(r-1)/2}\right),
&\mbox{if $r\geq 1$ is odd,}\\
\frac{(p+1)^2}{p^3-1}(p^{2r+1}+p^{2r-2}-2p^{r/2-1}\frac{1+p^2}{1+p}),
&\mbox{if $r\geq 2$ is even,}
\end{cases}
\EAL
for odd primes $p$. The corresponding Dirichlet series reads
\BAL
\Phi_{D_4}(s)&=\sum_{n=1}^\infty \frac{f_{D_4}^{}(n)}{n^s}
=\prod_{p\ne 2}\frac{1+p^{-s}+2p^{1-s}+2p^{-2s}+p^{1-2s}+p^{1-3s}}%
{(1-p^{2-s})(1-p^{1-2s})}\nn\\
&=1+\frac{16}{3^s}+\frac{36}{5^s}+\frac{64}{7^s}+\frac{152}{9^s}
+\frac{144}{11^s}+\frac{196}{13^s}+\frac{576}{15^s}+\frac{324}{17^s}+\cdots.
\nn
\EAL
Differences to Eq.~(\ref{frotpcub}) occur for all integers that are divisible
by the square of an odd prime.

\subsection{Primitive hypercubic lattice}

It is well known that the primitive hypercubic lattice $\Z^4$, which is a
sublattice of $D_4^*$ of index $2$, has a smaller point group than $D_4^*$,
containing only $192$ rotations, so that $[Aut(D_4^*):Aut(\Z^4)]=3$.
As a consequence, every class of symmetry
related rotations splits into three classes, one of which has the same
coincidence index as before,
$\varSigma_{\Z^4}(R)=\varSigma_{D_4}(R)$, while the other two classes have
index $2\varSigma_{D_4}(R)$. The number of coincidence rotations is thus
given by $192 f_{\Z^4}^{rot}(n)$, where $f_{\Z^4}^{rot}(n)$
is again multiplicative, but
slightly more complicated than $f_{D_4}^{rot}(n)$ (see~\cite{baa97,pzcsl2}).
It has the generating function
\BAL
\Phi_{\Z^4}^{rot}(s)&=\sum_{n=1}^\infty \frac{f_{\Z^4}^{rot}(n)}{n^s}
=(1+2^{1-s})\Phi_{D_4}^{rot}(s)
=(1+2^{1-s})\prod_{p\ne 2}\frac{(1+p^{-s})(1+p^{1-s})}{(1-p^{1-s})(1-p^{2-s})}
\nn\\
&=1+\frac{2}{2^s}+\frac{16}{3^s}+\frac{36}{5^s}+\frac{32}{6^s}
+\frac{64}{7^s}+\frac{168}{9^s}+\frac{72}{10^s}
+\frac{144}{11^s}+\frac{196}{13^s}+\frac{128}{14^s}
+\frac{576}{15^s}+\frac{324}{17^s}
+\cdots.
\label{frotccub}
\EAL
The CSLs split into two classes only, one of which has odd and the other
even index. Hence, the generating function for the number of CSLs reads
\BAL
\Phi_{\Z^4}(s)&=(1+2^{-s})\Phi_{D_4}(s)
=(1+2^{-s})\prod_{p\ne 2}\frac{1+p^{-s}+2p^{1-s}+2p^{-2s}+p^{1-2s}+p^{1-3s}}%
{(1-p^{2-s})(1-p^{1-2s})}\nn\\
&=1+\frac{1}{2^s}+\frac{16}{3^s}+\frac{36}{5^s}+\frac{16}{6^s}
+\frac{64}{7^s}+\frac{152}{9^s}+\frac{36}{10^s}
+\frac{144}{11^s}+\frac{196}{13^s}+\frac{64}{14^s}
+\frac{576}{15^s}+\frac{324}{17^s}
+\cdots.
\nn
\EAL
Differences to Eq.~(\ref{frotccub}) occur for all even integers and all
integers that are divisible by the square of an odd prime.
The spectrum is
$\sigma=\{2^\ell m\mid \ell\in\{0,1\}, \mbox{and $m\in\N$ odd} \}$.

\section{Root lattice $A_4$}

Usually, the $A_4$ lattice --- we will denote it by $L$ in the following ---
is embedded in $\R^5$ as a lattice plane.
However, this is inconvenient for our purposes and
we prefer to look at it in $\R^4$, since we want to exploit the useful
parameterisation by quaternions,
which we do not have at hand in $5$ dimensions.
A possible basis for $L$ consists of the 4 vectors
\BAL
(1,0,0,0), \frac{1}{2}(-1,1,1,1), (0,-1,0,0),
\frac{1}{2}(0,1,\tau-1,-\tau),
\EAL
where $\tau=\frac{1+\sqrt{5}}{2}$ is the golden mean whose algebraic conjugate
$\tau'$ can be written as $\tau'=-\frac{1}{\tau}=1-\tau$. 
The lattice $L$ cannot be identified with
a ring of quaternions. However, if we interpret the basis vectors as
quaternions, they relate to the icosian ring $\I$,
which is the $\Z[\tau]$--span
of the 4 quaternions
\BAL
(1,0,0,0),(0,1,0,0),\frac{1}{2}(1,1,1,1),\frac{1}{2}(1-\tau,\tau,0,1).
\EAL
Note that neither $L$ nor $\I$ are invariant under algebraic conjugation.
Combining the algebraic conjugation with a permutation
of the last two components yields an involution of the second kind
$\tilde{x}:=(x'_0, x'_1, x'_3, x'_2)$,
which was called twist map in~\cite{baaheu1,pzcsl4}. Note that $L=\tilde L$
is invariant under the twist map, which in addition is an 
antiautomorphism of $\I$. The twist map is the key to our analysis
since it provides us with a useful representation of $L$,
\BAL
L=\{x\in \I\mid x=\tilde{x}\}=\{x+\tilde{x}\mid x\in\I\}
\EAL
and its CSLs --- see below.

Any coincidence rotation of $L$ (and only those)
can be parameterised by a single primitive admissible 
quaternion $q\in\I$ by means of
$R(q)=\frac{1}{|q\tilde{q}|}qx\tilde{q}$ (see~\cite{pzcsl4}).
Here, admissible means that
$|q\tilde{q}|\in\N$.

Admissibility guarantees that one can define the factors
\BAL 
\alpha_q:=\sqrt{\frac{|\tilde q|^2}{\gcd(|q|^2, |\tilde q|^2)}}&&\mbox{and}&&
\alpha_{\tilde q}:=\alpha'_q=\sqrt{\frac{|q|^2}{\gcd(|q|^2, |\tilde q|^2)}},
\EAL
which are elements of $\Z[\tau]$. Note, however, that they are only
defined up to units of $\Z[\tau]$. In addition, we define the extension of
a primitive admissible quaternion by $q_\alpha=\alpha_q q$. The 
coincidence index now reads~\cite{pzcsl4}
\BAL
\varSigma_{A_4}(R(q))=
\frac{|q\tilde{q}|^2}{\gcd(|q|^2, |\tilde{q}|^2)}
=|q\tilde{q}|\alpha_q\alpha_{\tilde{q}}
=\lcm(|q|^2, |\tilde{q}|^2).
\EAL
Consequently, the coincidence spectrum is $\N$.
An explicit expression for the CSLs exists, too:
\BAL
L\cap\frac{q_\alpha L\tilde{q}_\alpha}%
{|q_\alpha\tilde{q}_\alpha|}
=\{q_\alpha x+\tilde{x}\tilde{q}_\alpha\mid
x\in\I\}
=(q_\alpha\I+\I \tilde{q}_\alpha)\cap L.
\EAL

Along the same lines as for the hypercubic lattices, one
can calculate the number
of coincidence rotations and the number of CSLs of a given index $n$. The
expressions obtained are slightly more complicated, since one has to
distinguish between the three cases $p=5$, $p=\pm 1\pmod 5$ and
$p=\pm 2\pmod 5$, as these primes behave differently and correspond
to the ramifying, splitting and inert primes of $\Z[\tau]$, respectively.
If $120f_{A_4}^{rot}(n)$ is the number of coincidence rotations ---
recall that the point group of $L$ contains $120$ elements --- its generating
Dirichlet series is given by
\BAL
\Phi_{A_4}^{rot}(s)&=\frac{1+5^{1-s}}{1-5^{2-s}}
\prod_{p\equiv\pm 1(5)}\frac{(1+p^{-s})(1+p^{1-s})}{(1-p^{1-s})(1-p^{2-s})}
\prod_{p\equiv\pm 2(5)}\frac{1+p^{-s}}{1-p^{2-s}}\nn\\
&=1+\frac{5}{2^s}+\frac{10}{3^s}+\frac{20}{4^s}+\frac{30}{5^s}
+\frac{50}{6^s}+\frac{50}{7^s}+\frac{80}{8^s}+\frac{90}{9^s}
+\frac{150}{10^s}+\frac{144}{11^s}+\cdots.\nn
\EAL
Similarly, we have for the number of CSLs
\BAL
\Phi_{A_4}(s)&=\left(1+6\frac{5^{-s}}{1-5^{2-s}}\right)
\prod_{p\equiv\pm 2(5)}\frac{1+p^{-s}}{1-p^{2-s}}
\prod_{p\equiv\pm 1(5)}\frac{1+p^{-s}+2p^{1-s}+2p^{-2s}+p^{1-2s}+p^{1-3s}}%
{(1-p^{2-s})(1-p^{1-2s})}
\nn\\
&=1+\frac{5}{2^s}+\frac{10}{3^s}+\frac{20}{4^s}+\frac{6}{5^s}
+\frac{50}{6^s}+\frac{50}{7^s}+\frac{80}{8^s}+\frac{90}{9^s}
+\frac{30}{10^s}+\frac{144}{11^s}+\cdots.\nn
\EAL
Here, a criterion analogous to Theorem~\ref{equalJ} exists, which will
be published elsewhere~\cite{pzcsl5}.

\section{Icosian ring}

It is also interesting to discuss the coincidences of the icosian ring $\I$
itself. Note that it is no lattice in 4--space
but a $\Z$--module of rank $8$.
However, this does not matter and one can argue similarly to the hypercubic
case. The main difference is that one has to work with $\Z[\tau]$
instead of $\Z$, which does not cause any problems. In fact, some steps
are even easier since the prime $2$ does not play a special role here.
So, a primitive quaternion is already reduced and hence the notion
of reducibility is not needed here.

Let us call a pair $(q,p)\in\I\times\I$ \emph{primitive admissible} if 
$q,p$ are primitive and $|qp|\in\Z[\tau]$. Then, a rotation is a coincidence
rotation if and only if it is parameterised by a primitive admissible pair
$(q,p)$. For any primitive admissible pair $(q,p)$, we can once more define
\BAL 
\alpha_q:=\sqrt{\frac{|p|^2}{\gcd(|q|^2, |p|^2)}}&&\mbox{and}&&
\alpha_{p}:=\sqrt{\frac{|q|^2}{\gcd(|q|^2, |p|^2)}},\nn
\EAL
which are again elements of $\Z[\tau]$, defined up to units. As before, 
$(q_\alpha,p_\alpha):=(\alpha_q q,\alpha_p p)$ is the corresponding
extension pair. We get the following expression for the coincidence
index~\cite{pzcsl5}:
\BAL
\varSigma_\I(R(q,p))
=\Nr(\lcm(|q|^2,|p|^2))
=\Nr(|q_\alpha|^2)
=\Nr(|p_\alpha|^2).\nn
\EAL
Here $\Nr(a)=|a a'|$ is the (number theoretic) norm of $a\in\Z[\tau]$.
Thus, the coincidence spectrum consists of all integers that contain
prime factors $p=\pm 2\pmod 5$ only with even power.
The CSMs read explicitly
\BAL
\I\cap \frac{q\I \bar{p}}{|qp|}
=q_\alpha\I + \I \bar{p}_\alpha.\nn
\EAL

If $7200 f_\I^{rot}(n)$ denotes the number of coincidence rotations
and $f_\I(n)^{}$ the number of CSMs, we get again $f_\I^{rot}(n)\geq f_\I(n)^{}$.
Again, equality fails to hold in general, and an analogue of
Theorem~\ref{equalJ} tells us which CSMs are equal~\cite{pzcsl5}.
The corresponding Dirichlet series read
\BAL
\Phi_{\I}^{rot}(s)&=\frac{(1+5^{-s})(1+5^{1-s})}{(1-5^{1-s})(1-5^{2-s})}
\prod_{p\equiv\pm 1(5)}
\left(\frac{(1+p^{-s})(1+p^{1-s})}{(1-p^{1-s})(1-p^{2-s})}\right)^2
\prod_{p\equiv\pm 2(5)}\frac{(1+p^{-2s})(1+p^{2-2s})}{(1-p^{2-2s})(1-p^{4-2s})}
\nn\\
&=1+\frac{25}{4^s}+\frac{36}{5^s}+\frac{100}{9^s}
+\frac{288}{11^s}+\frac{440}{16^s}+\frac{400}{19^s}
+\frac{900}{20^s}+\frac{960}{25^s}+\cdots
\EAL
and
\BAL
\Phi_{\I}(s)&=
\frac{1+5^{-s}+2\cdot5^{1-s}+2\cdot 5^{-2s}+5^{1-2s}+5^{1-3s}}%
{(1-5^{2-s})(1-5^{1-2s})}
\prod_{p\equiv\pm 1(5)}\frac{(1+p^{-s}+2p^{1-s}+2p^{-2s}+p^{1-2s}+p^{1-3s})^2}%
{(1-p^{2-s})(1-p^{1-2s})}\nn\\
&\times
\prod_{p\equiv\pm 2(5)}\frac{1+p^{-2s}+2p^{2-2s}+2p^{-4s}+p^{2-4s}+p^{2-6s}}%
{(1-p^{4-2s})(1-p^{2-4s})}\nn\\
&=1+\frac{25}{4^s}+\frac{36}{5^s}+\frac{100}{9^s}
+\frac{288}{11^s}+\frac{410}{16^s}+\frac{400}{19^s}
+\frac{900}{20^s}+\frac{912}{25^s}+\cdots.
\EAL

\section*{Acknowledgements}

The authors are grateful to Uwe Grimm, Manuela Heuer and Robert V. Moody
for helpful discussions on the present subject.  This work 
was supported by the German Research Council (DFG), within the CRC~701.




\label{lastpage}

\end{document}